\documentclass[12pt]{article}
\usepackage{amssymb,amsmath,amsthm,fancyhdr,array,calc}
\DeclareMathOperator{\GCD}{GCD}
\normalsize
\headwidth=\textwidth
\newcommand{\fancytitle}{\maketitle\thispagestyle{fancy}}

\title{\Large\bfseries On existence of Budaghyan-Carlet APN hexanomials}
\author{Antonia W. Bluher}
\date{August 2012}
\def\F{{\mathbb F}}
\def\Ker{{\rm Ker}}

\let\x=\times
\begin{document}
\fancytitle
\begin{abstract}
Budaghyan and Carlet \cite{BC} constructed a family of 
almost perfect nonlinear (APN)
hexanomials over a field with $r^2$ elements, and with terms of
degrees $r+1$, $s+1$, $rs+1$, $rs+r$, $rs+s$, and $r+s$, where $r=2^m$
and $s=2^n$ with $\GCD(m,n)=1$.
The construction requires a certain technical condition, which  
was verified empirically in a finite number of examples.
Bracken, Tan, and Tan \cite{BTT} proved the condition holds when $m\equiv 2$
or~$4\pmod 6$.
In this article, we prove that the construction of Budaghyan and Carlet
produces APN polynomials for all values of $m$ and $n$.

More generally, if $\GCD(m,n)=k\ge 1$, Budaghyan and Carlet showed that
the nonzero derivatives of the hexanomials are $2^k$-to-one maps from 
$\F_{r^2}$ to $\F_{r^2}$, 
provided the same technical condition holds. We prove their construction 
produces polynomials with this property for all $m$ and~$n$. 
\end{abstract}

\section{Introduction}

If $f$ is a function from $\F_{2^n}$ to $\F_{2^n}$, 
one can ask for the number of
solutions to $f(x+a) = f(x) + b$, where $a,b\in\F_{2^n}$ and $a$ is nonzero.  
Note that if $x$ is one solution, then $x+a$ is another, so the solutions
come in pairs. The function $f$ is said to be {\it almost perfect nonlinear}
(APN) if there are always exactly zero or two solutions. 
The function $f(x+a) + f(x)$ is called a {\it derivative} of $f$.
An APN function is simply a function whose derivatives yield two-to-one
maps on $\F_{2^n}$. As pointed out by Carlet, Charpin, and Zinoviev \cite{CCZ},
the APN property is equivalent to the property that a certain binary linear 
code defined in terms
of $f$ is double error-correcting. Construction of APN functions is
a recurring theme in the literature; see \cite{BDKM}, \cite{BDMW}, \cite{CarletDCC},
and the survey article \cite{Carlet}.

Let $r=2^m$ and $s=2^n$, where $m,n \ge 1$. 
For $d\in \F_{r^2}\setminus \F_r$, and $c\in\F_{r^2}$,
Budaghyan and Carlet \cite{BC} consider the hexanomial
\begin{equation}
F(x) = x(x^s + x^r + c x^{rs}) + x^s(c^r x^r + d x^{rs}) + x^{(s+1)r}. \label{Fdef}
\end{equation}
For any positive integer
$N$, denote by $\mu_N$ the group of $N$-th roots of unity in the algebraic
closure of $\F_2$.  If $M$ is
odd, then $\mu_M$ has order~$M$, and $\mu_M\subset \mu_N$ if and only if 
$M$ divides $N$. In particular, $\mu_{r+1}\subset \mu_{r^2-1} = \F_{r^2}^\x$,
where $\F^\x$ denotes the nonzero elements of a field $\F$.

\outer\def\proclaim#1. #2 \par{\medbreak\noindent{\bf #1.\enspace}{\it #2}\par
\ifdim\lastskip<\medskipamount\removelastskip\penalty55\medskip\fi}
\newenvironment{pf}{\noindent {\bf Proof.}}{\qed\vskip 6pt}

\proclaim{Theorem 1 (Budaghyan and Carlet \cite{BC})}.
{
If $y^{s+1} + c y^s + c^r y + 1$ has no roots $y$ belonging to $\mu_{r+1}$,
then all the derivatives of $F$
are $2^k$-to-1 mappings from $\F_{r^2}$ to $\F_{r^2}$, where  $k=\GCD(m,n)$. 
(In particular, if $k=1$, then $F(x)$ is APN.)
}

Let us say that the pair $(r,s)$ is {\it BC-compatible} if $c \in \F_{r^2}$
can be found satisfying the hypothesis of the theorem.
We found an exact and surprisingly simple criterion for BC-compatibility:

\proclaim{Theorem 2}.  {The pair $r=2^m$, $s=2^n$ is BC-compatible if and only if
$m > 1$ and $n/m$ is not an odd integer. 
}

Previously, it was known that $(r,s)$ is BC-compatible only in some
special cases.
In \cite{BC} it was found empirically that $(2^m,2)$ is BC-compatible
whenever $6\le 2m\le 500$ and $3 \nmid m$, and also in at least
140 of the 166 checked cases when $3$ divides $m$.  Later,
Bracken, Tan, and Tan \cite{BTT} proved that $(r,s)$ is 
BC-compatible if $m \equiv 2$ or $4\pmod 6$ and $\GCD(m,n)=1$, and in particular the
Budaghyan-Carlet APN hexanomials belong to an infinite family.
All the cases in \cite{BC} and \cite{BTT} 
satisfy that
$y^{s+1} + c y^s + c^r y + 1$ has no roots in $\F_{r^2}$.
This condition is stronger than the required hypothesis, 
since $\mu_{r+1} \subset \F_{r^2}$.

Theorem~2 implies that if $c$ is properly selected, then $F(x)$ is APN whenever
$m>1$ and $\GCD(m,n)=1$. We will show that $F(x)$ is APN when $m=1$ also, so 
in fact the only requirement is $\GCD(m,n)=1$.  More generally, we prove the following.

\proclaim{Theorem 3}. {For all $r = 2^m$ and $s=2^n$, and for all 
$d\in\F_{r^2}\setminus\F_r$, a value $c\in\F_{r^2}$ can be found such that 
all the nonzero derivatives of $F(x)$ are $2^k$-to-one mappings from
$\F_{r^2}$ to $\F_{r^2}$, where $k=\GCD(m,n)$.}

For another viewpoint on the APN hexanomials $F(x)$, see 
\cite[Section 4.2.1]{CarletDCC}, where it is shown that they belong to a 
family that is constructed using bent functions.

\section{Proof of Theorem 1}

For completeness, we present the proof by Budaghyan and Carlet of Theorem~1.
As above, $r=2^m$, $s=2^n$, $d \in \F_{r^2}\setminus \F_r$, and
$c \in \F_{r^2}$.  Note that $\F_r\cap\F_s=\F_u$, where $u=2^k$, $k=\GCD(m,n)$.
Let $F(x)$ be the hexanomial defined in (\ref{Fdef}).
Assuming the hypothesis that $y^{s+1} + c y^s + c^r y + 1$ has no
roots in $\mu_{r+1}$, we are to show that for any nonzero $a \in \F_{r^2}$
and any $b \in \F_{r^2}$, the equation
$$F(x) + F(x+a) = b$$
has exactly zero solutions or exactly $u$ solutions in $\F_{r^2}$.

Denote the
number of solutions by $N(a,b)$. Let $G_a(x)=F(ax) + F(ax + a) + F(a).$
Then $N(a,b)$ is the number of solutions in $\F_{r^2}$ to $G_a(x)=F(a)+b$.
We claim that $G_a$ is an $\F_u$-linear function.  
Accepting this for the moment,
we see that proving $N(a,b) \in \{0,u\}$ is equivalent to showing
that $\Ker(G_a)$ (considered as an $\F_u$-linear function on $\F_{r^2}$)
has order $u$. We will in fact show $\Ker(G_a) = \F_u$.

To see that $G_a$ is $\F_u$-linear, we note that the terms in $F(ax)$ are
of the form $\alpha x^{v+w}$ or $\alpha x^v$, where $\alpha \in \F_{r^2}$
and $v,w\in \{r,s,rs,1\}$ (all powers of $u$).  Thus, 
$G_a$ is a sum of terms $\alpha ( x^{v+w} + (x+1)^{v+w} + 1) = \alpha(x^v + x^w)$.
This is $\F_u$-linear because $v$ and $w$ are powers of $u$. Note also that
$\Ker(G_a)$ contains $\F_u$, because $x^v+x^w=x+x=0$ for all $x\in\F_u$.

Now $G_a(x) = a^{s+1}(x+x^s) + a^{r+1}(x+x^r) + ca^{rs+1} (x+x^{rs})
+ c^r a^{r+s} (x^r+x^s) + d a^{s+rs} (x^s+x^{rs}) + a^{(s+1)r}(x^{rs}+x^r)$.
Suppose $G_a(x_0)=0$ with $x_0\in \F_{r^2}$.   Then of course
$G_a(x_0) + G_a(x_0)^r = 0$. Using that $x_0^{r^2}=x_0$,
$a^{r^2}=a$, $c^{r^2}=c$, $d^{r^2}=d$, we find that many terms in $G_a(x_0)^r$
cancel with terms in $G_a(x_0)$. The result is
$$0=G_a(x_0)+ G_a(x_0)^r = (d+d^r) a^{s+rs} (x_0+x_0^r)^s.$$
Now $d+d^r\ne0$ since $d\not\in \F_r$, $a^{s+rs}\ne0$ since $a\ne 0$.
So we have $x_0 + x_0^r=0$.  Returning to the original formula for $G_a$ 
and using the relation $x_0 = x_0^r$,
we see that every term either vanishes or becomes a multiple of $x_0+x_0^s$:
\begin{eqnarray*}
0&=& G_a(x_0)\\
&=&  (x_0+x_0^s) (a^{s+1} + c a^{rs+1}  + c^r a^{r+s} + a^{(s+1)r} ) \\
&=& (x_0 + x_0^s) a^{s+1} ( 1 + c a^{(r-1) s} + c^r a^{r-1} + a^{(s+1)(r-1)} ).
\end{eqnarray*}
Since $a$ is nonzero, the term $a^{s+1}$ is nonzero. Since $a^{r-1}$ belongs
to $\mu_{r+1}$, the hypothesis of the theorem implies that
$1 + c a^{(r-1) s} + c^r a^{r-1} + a^{(s+1)(r-1)} $ is nonzero.
So we conclude that $G_a(x_0) = 0$ implies $x_0^r=x_0$ and $x_0^s=x_0$, 
{\it i.e.} $x_0 \in \F_r \cap \F_s = \F_u$. This proves that $\Ker(G_a)=\F_u$,
as claimed.

\section{Proof of Theorem 2}

As above, let $r=2^m$ and $s=2^n$, where $m,n \ge 1$. Let
$$G(c,y)=y^{s+1} + c y^s + c^r y + 1.$$
The technical condition needed in Theorem~1 for the hexanomial $F(x)$ 
to have desired properties is that there exists $c \in \F_{r^2}$ such that 
$G(c,y)$ has no roots in $\mu_{r+1}$.  If such $c$ exists, then we say
that the pair $(r,s)$ is BC-compatible.  We first need a lemma.

\proclaim{Lemma 1}. {$r+1$ divides $s+1$ if and only if $n/m$ is an odd integer.
}

\begin{proof} 
First, suppose $n/m=\ell$ is an odd integer, and we will show
that $r+1$ divides $s+1$. Since $\F_{2^a}\subset \F_{2^b}$ if and
only if $a|b$, and since $2m|2n$, we see that $\F_{r^2}\subset \F_{s^2}$. Since
$x\in \F_{2^a}^\x$ if and only if the order of $x$ divides $2^a-1$, we see that
$\mu_{r+1} \subset \F_{r^2}$ and $\mu_{s+1} \subset \F_{s^2}$. Let $\tau$ denote
the Frobenius map on $\F_{s^2}$ (given by squaring), 
$\rho=\tau^m$, and $\sigma = \tau^n = \tau^{m\ell}  = \rho^\ell$. 
Note that $\rho(a)=a^r$ and $\sigma(a)=a^s$, for $a\in \F_{s^2}$. Now
\begin{equation}
\mu_{r+1} = \{z\in \F_{s^2}^\x : \rho(z) = 1/z \},
\qquad \mu_{s+1} = \{z\in \F_{s^2}^\x : \sigma(z) = 1/z \}.   \label{FrobFormula}
\end{equation}
Since $\ell$ is odd, we see that if $z\in \mu_{r+1}$ then $\sigma(z)
= \rho^\ell(z)=1/z$, and so $z\in \mu_{s+1}$.  Thus, $\mu_{r+1} \subset
\mu_{s+1}$, and consequently $r+1$ divides $s+1$.

To prove the converse, suppose that $r+1$ divides $s+1$ and we will prove
that $n$ is an odd multiple of~$m$. Let $K_r$ denote
the subfield of the algebraic closure of $\F_2$ that is generated 
by $\mu_{r+1}$. We claim
$K_r = \F_{r^2}$.  First, $\mu_{r+1} \subset \mu_{r^2-1} = \F_{r^2}^\x$, so
$K_r \subset \F_{r^2}$. Now $\F_{r^2}$ can be viewed as a vector space
over $K_r$. If the dimension is $d$, then $r^2=|K_r|^d \ge (r+1)^d > r^d$. 
So $d=1$, and consequently $K_r=\F_{r^2}$ as claimed.

Since $r+1$ divides $s+1$, we have $\mu_{r+1} \subset \mu_{s+1}$, so the field
generated by $\mu_{r+1}$ is contained in the field generated by $\mu_{s+1}$.
That is, $\F_{r^2}=\F_{2^{2m}} \subset \F_{s^2} = \F_{2^{2n}}$. It follows
that $m$ divides $n$, say $n = \ell m$.  Let $\tau$, $\rho$, $\sigma$ be
as above, and let $1\ne z \in \mu_{r+1}$. By (\ref{FrobFormula}),
$\rho(z)=1/z$.  Since $\sigma = \rho^\ell$, and $z\ne 1/z$,
 we see that $\sigma(z)=1/z$ if $\ell$ is odd,
and $\sigma(z) = z\ne 1/z$ if $\ell$ is even.  On the other hand, 
$z\in\mu_{r+1}\subset \mu_{s+1}$, so by (\ref{FrobFormula}),  
$\sigma(z)=1/z$.  Then $\ell$ must be odd.
\end{proof}

Now we prove our theorem.

\proclaim{Theorem 2}.  {Let $r$ and $s$ be arbitrary positive integral powers of two,
and let 
$$G(c,y) = y^{s+1} + c y^s + c^r y + 1.$$  
There exists $c\in\F_{r^2}$ such that  $G(c,y)$ has no roots in 
$\mu_{r+1}$ if and only if $r>2$ and $r+1$ does not divide $s+1$.
(By the lemma, these conditions on $r$ and $s$ are equivalent to 
$m>1$ and $n/m$ is not an odd integer.)
}

\begin{proof}
First let us show if $r = 2$ then $G(c,y)$ has a root in $\mu_3$
for any $c\in \F_4$.
If $c\in \{0,1\}$ then $G(c,1)=0$.
If $c\in\F_4\setminus \F_2$ then $G(c,y)=0$  for $y=c \in \mu_3$.
This establishes the result when $r=2$.

Now let us show if $r + 1$ divides $s + 1$ then for all $c\in\F_{r^2}$, the
polynomial $G(c,y)$ has a root $y\in \mu_{r+1}$.
If $c =0$, then $G(c,1)=0$.
If $c\ne 0$, then set $y=c^{(r/2)(r-1)}$.
This belongs to $\mu_{r+1}$, because
$y^{r+1}=(c^{r/2})^{r^2-1}=1$.
Since $r+1$ divides $s+1$, we have $y^{s+1}=1$, so 
$$ G(c,y) = 1 + c/y + c^r y + 1 = (c/y)(1+c^{r-1}y^2) =  (c/y)(1+c^{r^2-1})=0.$$

For the remainder of the proof, assume $r>2$ and $r+1\nmid s+1$. We must find
$c \in \F_{r^2}$ such that $G(c,y)$ has no roots $y\in\mu_{r+1}$.
For $y \in \mu_{r+1}$, let
$$X_y  = \{ a \in \F_{r^2} : G(a,y) = 0 \}.$$ 
We are seeking $c \in \F_{r^2} \setminus X$, where
$$X = \cup_{y \in \mu_{r+1}} X_y.$$
Such $c$ exists if and only if $|X| < r^2$.

Since $G(c,y)$ has degree $r$ in the variable $c$, we have $|X_y| \le r$.    
This gives a bound: 
$$|X| \le \sum_{y\in \mu_{r+1}} |X_y| \le r(r+1).$$
This bound is not good enough, as we need to show $|X| < r^2$.
To attain this, we must take into account that the sets $X_y$ are not disjoint.

We consider separately the two cases: $r+1$ divides $s-1$, and $r+1$ does
not divide $s-1$.  If $r+1$ divides $s-1$, then for $y\in \mu_{r+1}$ we have
$G(c,y)=y^2 + (c  + c^r) y + 1$.
It follows that
$$G(c,1/y) = y^{-2} + (c+c^r) y^{-1} + 1 = y^{-2} G(c,y),$$
and so $X_y = X_{y^{-1}}$. Consequently,
$X =  \cup X_y$, where the union includes
one representative among each pair $\{y,1/y\}$.
There are $1 + r/2$ representatives,  giving $|X|\le r (1 + r/2)$. 
Since $r>2$ by hypothesis, this is less than $r^2$, as required.

Finally,  we consider the case where $r+1 \nmid s+1$ and
$r+1 \nmid s-1$. Observe that $X_1 = \{ a \in \F_{r^2} : 1 + a + a^r + 1 = 0\}
= \F_r$.  Also, observe that if $y \in \mu_{r+1}$ then $G(y,y)=0$,
so $y \in X_y$.  Thus, $X_1 \subset Z \subset X$, where
$$Z = \F_r \cup \mu_{r+1}.$$  
It follows that
$$X = Z \cup  \left( \cup_{y\in\mu_{r+1}, y \ne 1} X_y \setminus Z\right),$$
and so
\begin{eqnarray*}
|X| &\le& |Z| + \sum_{y\in \mu_{r+1}, y \ne 1} |X_y \setminus Z| \\
    &=&  2r + \sum_{y \in \mu_{r+1}, y \ne 1} \left( |X_y| - |X_y \cap Z| \right) \\
    &\le& 2r + \sum_{y\in \mu_{r+1}, y\ne 1} \left( r - |X_y \cap Z| \right) \\
    &=& 2r + r^2 - \sum_{y \in \mu_{r+1}, y \ne 1} |X_y \cap Z|.
\end{eqnarray*}   
This leads to the inequality
\begin{equation}
r^2 - |X| \ge  \sum_{y \in \mu_{r+1}, y \ne 1} \left(|X_y \cap Z|-2\right).
\end{equation}
So to demonstrate that $|X|<r^2$,  it suffices
to show that $|X_y \cap Z| \ge 2$ for all $y \in \mu_{r+1}\setminus\{1\}$,
and $|X_y \cap Z| > 2$ for at least one $y$.
We will do this by constructing some explicit elements of $X_y \cap Z$.

Two elements of $X_y \cap Z$ are $y$ and $y^{-s}$. These are in $X_y$ because
for $c = y$, 
$$y^{s+1} + c y^s + c^r y + 1 = y^{s+1} + y^{s+1} + y^{r+1} + 1 = 0,$$
and for $c=y^{-s}$, 
$$y^{s+1} + c y^s + c^r y + 1 = y^{s+1} + 1 + y^{s+1} + 1 = 0.$$
Note that $y$ and $y^{-s}$ are distinct if and only if $y^{s+1}\ne 1$.

If $y^{s-1} \ne 1$ then we can obtain another element of $X_y \cap Z$ by setting
$$c_0 = (y^{s+1} + 1)/(y^s + y).$$
Here $c_0 \in \F_r$, because (using $y^r = 1/y$) we have
$$c_0^r = (y^{-(s+1)} + 1)/(y^{-s}+y^{-1}) = (1 + y^{s+1})/(y + y^s) = c_0.$$
Also $c_0 \in X_y$, because 
$$y^{s+1} + c_0 y^s + c_0^r y + 1 =  (y^{s+1} + 1) + c_0 (y^s + y) = 0.$$
Since $c_0 \in \F_r$ and $\F_r \cap \mu_{r+1} = \{1\}$, we know
$c_0$ is distinct from $y$ and $y^{-s}$.  

In summary, for $y \in \mu_{r+1} \setminus \{1\}$ we have:
\begin{itemize}
\item If $y^{s-1} \ne 1$ and $y^{s+1} \ne 1$, then $c_0$, $y$, and $y^{-s}$ 
are distinct elements of $X_y\cap Z$. 
\item If $y^{s-1} \ne 1$ but $y^{s+1} = 1$, then $c_0$ and $y$ are distinct
elements of $X_y\cap Z$.
\item If $y^{s-1} = 1$ then $y$ and $y^{-s}$ are distinct elements of $X_y\cap Z$.
\end{itemize}

We see that $|X_y \cap Z| \ge 2$ always.  Moreover, when $y$ is a primitive
$(r+1)$th root of unity, then from the hypothesis that $r+1$ does not divide
$s+1$ or $s-1$, we will have that $y^{s+1} \ne 1$ and $y^{s-1} \ne 1$, so
$|X_y \cap Z| \ge 3$.  As noted above, this  completes the demonstration
that $|X|<r^2$, and completes the proof. 
\end{proof}

\section{Proof of Theorem 3}

Theorem 3 asserts that for $r=2^m$ and $s = 2^n$, and any choice of
$d\in\F_{r^2}\setminus\F_r$, there always exists $c\in\F_{r^2}$ such that
the nonzero derivatives of the hexanomial $F(x)$ given by (\ref{Fdef}) are 
$2^k$-to-one mappings from $\F_{r^2}$ to $\F_{r^2}$, where $k=\GCD(m,n)$. 
Here we provide a proof.

If $m$ does not divide $n$, then $(r,s)$ is BC-compatible by Theorem~2, so 
Theorem~3 holds.  If $m$ divides~$n$, then the next lemma shows that any
choice of $c$ will work, so that Theorem~3 again holds.

\proclaim{Lemma 2}. {If $m$ divides $n$ (equivalently, $\F_r\subset\F_s$), 
then the nonzero derivatives of $F(x)$ are $r$-to-one mappings from
$\F_{r^2}$ to $\F_{r^2}$, for any choice of $c\in\F_{r^2}$ and 
$d\in\F_{r^2}\setminus\F_r$.}

\begin{proof} For nonzero $a\in\F_{r^2}$, let $G_a(x) = F(ax) + F(ax+a) + F(a)$.
As explained in the proof of Theorem~1, it suffices to prove that $G_a$ has exactly
$r$ roots in $\F_{r^2}$. If $x_0\in\F_{r^2}\setminus\F_r$, then using the
relation $x_0^{r^2}=x_0$, we find that $G_a(x_0)+G_a(x_0)^r = (d+d^r)a^{s+rs}
(x_0+x_0^r)^s$. This is nonzero, therefore $G_a(x_0)\ne0$. If $x_0\in\F_r$, then 
using the relation $x_0^r=x_0$ we find that $G_a(x_0) = (x_0+x_0^s) a^{s+1} 
(1 + c a^{(r-1)s}+ c^r a^{r-1} + a^{(s+1)(r-1)})$. Since $x_0\in\F_r\subset\F_s$,
we see that $x_0+x_0^s=0$, and so $G_a(x_0)=0$. This establishes that $G_a$ has 
exactly $r$ roots in $\F_{r^2}$, as required.
\end{proof}

\end{document}